\documentclass[12pt]{amsart}
\usepackage{amsmath,amsthm,amsfonts,amssymb}
\usepackage{color}
\usepackage{amsfonts,amsthm,amsmath,amssymb,verbatim}
\usepackage{graphicx}
\usepackage[active]{srcltx}
\usepackage[all,cmtip]{xy}
\usepackage{tikz}
\usepackage{xypic}
\usepackage{hyperref}
\hypersetup{
    colorlinks=true,
   linkcolor=blue,
   filecolor=magenta,
    urlcolor=cyan,
}

\newtheorem{thm}{Theorem}[section]
\newtheorem{prop}[thm]{Proposition}
\newtheorem{lemma}[thm]{Lemma}

\newtheorem{cor}[thm]{Corollary}

\theoremstyle{remark}

\newtheorem{defin}{Definition}



\def\C{\mathbb{C}}

\def\Q{\mathbb{Q}}
\def\Z{\mathbb{Z}}

\def\F{\mathbb{F}}

\def\lra{\longrightarrow}

\def\eps{\varepsilon}

\usepackage{graphicx}  

\title{Generalized Sato-Tate   and  quadratic residues}
\author{Sergey Vl\u adu\c t  }\thanks{Aix--Marseille Universit\'{e}, CNRS, I2M UMR 7373, Marseille, France, sergevladuts@ya.ru, vladut@iml.univ-mrs.fr}\address{Aix--Marseille Universit\'{e}, CNRS, I2M UMR 7373, Marseille, France, sergevladuts@ya.ru, vladut@iml.univ-mrs.fr}
\date{}

\begin{document}

\maketitle 

\section{Introduction}
The present note is a continuation of \cite{KTVZ} and deals with several questions concerning the distribution of consecutive quadratic residues that were left open there.

One of the main results of \cite{KTVZ} provides rather precise information on the distribution of consecutive quadratic residues modulo very large primes. These results are based on strong statistical independence assumptions that appear natural but are, at present, completely out of reach.

The main observation in this note is that these assumptions follow from certain particular cases of the Generalized Sato–Tate conjecture (GST). Moreover, some already established cases of the GST yield interesting upper and lower bounds for the corresponding arithmetic functions.

Section 2 recalls the formulation of the GST.
Section 3 summarizes several results from \cite{KTVZ} concerning consecutive quadratic residues and their connections with algebraic curves and abelian varieties over finite fields.
Section 4 presents asymptotic results on the distribution of four and five consecutive quadratic residues, under the assumption of the GST for certain products of elliptic curves of dimension at most five.
Finally, Section 5 establishes unconditional bounds, similar in form but weaker than those of Section 4, by using known cases of the GST for products of two elliptic curves.

\smallskip\smallskip {\bf Acknowledgements.} I would like to thank Andrew Sutherland for his very useful remark \cite{Su1}. Many thanks to David Kohel for his help with numerous calculations, especially concerning the curve $C$ on page 8.
 
\section{The Generalized Sato-Tate conjecture}

Let us recall the   setting for the (generalized)
Sato-Tate conjecture and some known elementary cases. Recall that  we consider only the   simplest case of the base field $\Q.$\smallskip

Thus, let $A$ be an abelian variety over $\Q$ of dimension $g\ge 1$. The Generalized Sato-Tate conjecture (GST, or GST(A)): \smallskip \smallskip

 I. associates with $A$ a (conjugacy class of) closed algebraic subgroup\linebreak $\mathrm {ST}_A<\mathrm{Sp}_{2g}(\C)$;\smallskip\smallskip

 II. states that the family $\{\mathrm {F}_{p},p\notin S_A\}\subset \mathrm{ ST}_A,$ where $p$ is a (rational) prime, formed by natural images of $p$-Frobenius  elements $\mathrm {Fr}_{p}\in \mathrm{ G}_\Q := \mathrm{ Gal}({\bar \Q}/\Q)$ in $\mathrm{ ST}_A$ is equidistributed  with respect to  the  probabilistic Haar measure $\mu_A$ on $\mathrm{ ST}_A,\;S_A$ being the finite set of bad reduction primes for $A.$   \smallskip\smallskip

 We recall then the construction of $\mathrm {F}_{p}$ and the notion of being equidistributed.
 
 Fix a prime $\ell\ne p$ and let 
 $$ \phi_\ell:\mathrm{ G}_\Q \longrightarrow \mathrm{ GL}_{2g}(\Q_\ell)\subset \mathrm{ GL}_{2g}(\bar\Q_\ell)$$
 be  the standard $\ell$-adic representation, given by the Galois action on torsion points of $A.$ Choose and fix an isomorphism $I:\bar\Q_\ell\lra \C, $ and hence  an isomorphism 
 $$I_{2g}:\mathrm{ GL}_{2g}(\bar\Q_\ell)\lra\mathrm{ GL}_{2g}(\C),
 $$
 then by Weil-Deligne we have 
 $$\mathrm {F}_{p}:=I_{2g}(\phi_\ell(\mathrm {Fr}_{p})/ \sqrt  p)\in \mathrm{Sp}_{2g}(\C)
 $$ 
 and, in fact, 
$  \{\mathrm {F}_{p}, p\notin S_A\}\subset\mathrm{ ST}_A$ (and  Zarisky-generate $\mathrm{ ST}_A$).\smallskip

Let then $T$ be a compact and $\mu:=\mu_T$ be a probabilistic Borel measure on $T.$ If $S=\{t_1, t_2, \ldots t_n, \ldots\}\subset T$ then $S$ is equidistributed  with respect to $\mu$, or $\mu$-equidistributed if the sequence $\{\lambda_n\}$ of discrete measures
$$ \lambda_n:=\frac{1}{n}\sum_{i=1}^n \delta_{i}$$
*-weakly converges to $\mu$ for $n$ going to infinity, $\delta_{i}$ being the Dirac measure supported on  $t_i.$
Recall that it means that for any continuous function
$f :T\lra \C$ one has 
$$ \lim_{n\to \infty} \frac{1}{n}\sum_{i=1}^n f(t_{i})=\int_T f d\mu\, . $$

Since one has the canonical surjection 
$$ \Phi:\mathrm{Sp}_{2g}(\C)\lra [-2g, 2g], M\mapsto  \mathrm{Tr}(M),$$
we can consider the projected measure $\nu_A:=\Phi_{*}(\mu_A)$ on $[-2g, 2g],$ and the sequence or traces  $\{\Phi(\mathrm{F}_p), p\notin S_A\} $ is $\nu_A$-equidistributed.\smallskip\smallskip

{\em Example 1.} Let $g=1$, i.e. $A=E,$ an elliptic curve over $\Q$ and $\mathrm{ ST}_E\subset \mathrm{Sp}_{2}(\C) = \mathrm{SU}_{2}(\C)  .$ Then there are two possibilities: \smallskip\smallskip

a). $E$  is a CM-curve, that is $\mathrm{End}_{\C}(E)\ne\Z;$ then $\mathrm{ ST}_E=N(T),$ the normalizer of the diagonal subgroup $T$ which has index 2 in $N(T).$
\smallskip\smallskip

b).  $\mathrm{End}_{\C}(E)=\Z$, that is, $E$ has no complex multiplication, and $\mathrm{ ST}_E=\mathrm{SU}_{2}(\C).$ \medskip
Note that in any case $\mathrm{End}_{\Q }(E)=\Z.$\smallskip

The case a) is well-known (Hecke,Weil,...) since sixties,
while the equidistribution in the case b) is just the initial Sato-Tate conjecture for elliptic curves, proved in \cite{HSBT}.

A  simple and well-known calculation (see, e.g., \cite{Se}, p.108) shows that 

\begin{equation}\label{nu}\nu_1= \nu_{CM} =\frac12 \delta_0+\frac{dt}{2\pi\sqrt{4-t^2}},\quad  \nu_2= \nu_{nCM} = \frac{\sqrt{4-t^2}dt}{2\pi}.\end{equation}

In the present note we need only  very special cases of
GST  when $A$ is a product of elliptic curves.

In particular, for $g=\dim A=2$,  we use  $A= E_1\times E_2$ 
for 2 non-isogenous EC $E_1,E_2;$ note that in this  case GST(A) is known  thanks to \cite{Jo}. \smallskip

For $g=3$  we need $A=E_0\times  E_1\times E_2,$ where $E_0$ has CM, and $E_1,E_2$, being non-isogenous (over $\C$) are without CM.  This  case is not described in the literature, but it could  be accessible \cite{Su1} using the methods of  \cite{Jo}. \smallskip

However, the next case of interest for us has $g=17$ and thus is  far from the current proof methods. \smallskip

Nevertheless, in the rest of our note we suppose,  except in
Section 5 that GST(A) is true for all  abelian varieties $A$ used below.

\section{Consecutive quadratic residues and algebraic curves}
Let $p$ be an odd prime.
Consider the sequence $1$, $2$, \ldots, $p-1$.
Replace every number $i$ in the sequence by the letter $R$ if $i$ is a quadratic residue modulo $p$, and by the letter $N$ otherwise.
Denote by $W_p$ the resulting word.
\begin{defin}
Let $S$ be a word of length $t\le{p-1}$ that contains only $R$ and $N$.
Define $n_p(S)$ as the number of sub-words of $W_p$ that coincide with $S$, and are formed by $t$ consecutive elements of $W_p$. The word $$\underbrace{R\ldots R}_{t}$$  will be also denoted by $R^t$.
\end{defin}
 
Below we will concentrate only on the behaviour of  $n_p(R^t)$ and denote it by $n_p(t)$.  For other cases and some detail see our  paper \cite{KTVZ}.\smallskip 
 
Note that if $t=1$, then $n_p(1)=\frac{p-1}2$,
 i.e., the number of quadratic residues is equal to the number of non-residues. 
 
Aldanov's paper \cite{A} gives the answer for $t=2:$ 
$$n_p(2)=\frac{p-5}{4}, \hbox{ if } p=4k+1,  \;n_p(2)= \frac{p-3}{4}, \hbox{ if } p=4k+3.$$
for all (odd) $p$. 
 
Jacobstahl's 1906 thesis \cite{J}  explores the case $t=3$, which is  trickier.
However, for $p=4k+3$, the answer is simple:
$$n_p(3) = \frac{p-3-2\genfrac(){0.5pt}{0}{2}{p}}{8}.$$

Here and below  $\genfrac(){0.5pt}{0}{a}{p}$ is the Legendre symbol.\smallskip 

On the contrary, for $p=4k+1$  the formula is not so simple.
Namely, define  $J(k)\in 2\Z$ by  
$$a:=a(p):=J(k)=\sum_{i=1}^{4k-2}\genfrac(){0.5pt}{0}{i(i+1)(i+2)}{p}.$$
 Then
$$n_p(3)= \frac18\left\{p-11-4\genfrac(){0.5pt}{0}{2}{p}\right\} +\frac{a}4.$$
 
  For an  arbitrary $t\ge 1,$ and an odd $ p $ one can prove \cite[p.3]{C}, see also \cite[section 2.3]{KTVZ}:   

\begin{equation}\label{con} 
   n_p(t)=2^{-t}\sum_{j=1}^{p-t-1}\prod_{i=1}^{t} \left(1+\genfrac(){0.5pt}{0}{i+j -1}{p}\right)   
\end{equation}
which generalizes   Jacobsthal's formula.

Then a simple calculation, see \cite[section 2.3]{KTVZ} gives:

$$ n_p(t)=2^{-t}p +2^{-t}\sum_{s=1}^{t} \sum_{\substack{T\subseteq[1,t]\\ |T|=s}}N_T +c_p $$ 
 for a certain constant $c_p=c_p(t)\in 2^{-t} \Z$ which we will mostly ignore. Here 
 $$N_T:=\sum_{j\in \F_p}\prod_{i\in T}\genfrac(){0.5pt}{0}{i+j -1}{p}=\sum_{j\in\F_p}\genfrac(){0.5pt}{0}{f_T(j)}{p}\in \Z,$$
 for
 $$f_T(X):=\prod_{(i+1)\in T}(X+i)\in\Z[X].$$

 If $C_T$ is a hyper-elliptic curve 
  $$C_T: y^2=f_T(x), \quad, $$
  of genus $g_T=g(C_T)=\lfloor(s-1)/2\rfloor$ for $s=|T|,$
then $p+1-N_T$ is just the trace of   Frobenius  for $C_T$. 

Note also that for $t\ge 2$ one has $$g_t:=\sum_{\substack{T\subseteq[1,t]\\ }}g_T=2^{t-2}(t-3)+1,$$ see \cite[section 2.3]{KTVZ} \medskip
 
Let us look now at small values of $t$.\medskip

$\bf{t=1} .$ Here there is only one $T$ with $ f_T=X,$
$C_T$ is rational, $N_T=0$ and $c_p=-\frac12.$ 
\smallskip

$\bf{t=2} .$ Here $g_2=0,$ there are 3 possible $T,$ all with rational $C_T,\,N_T=0$ and
$c_p=-\frac14$ or $c_p=-\frac34$ depending on $p\mod 4.$\medskip

 $\bf{t=3} .$ Here $g_3=1;$ there are 7 possible $T,$ all of them, except for $\tilde T=\{1,2,3\}$ having rational $C_T,\,N_T=0,T\ne \tilde T;$ while $ C_{\tilde T}$ is the  "Gauss" CM elliptic curve $E_0$
$$y^2=x(x+1)(x+2),\; y^2=z^3-z,\; z:=x+1$$
which leads to  Jacobsthal's formula with $c_p, |c_p|\le \frac32$ depending on $p\mod 8.$\medskip

{\bf $t$=4.} Here $g_4=5$ and we get 5 non-rational curves,  all elliptic (4 for $|T|=3$ and 1 for $|T|=4$):
 \begin{equation}\label{ec}E_0: y^2=x(x+1)(x+2); \quad E_1: \; y^2=x(x+1)(x+3); \quad E_2: \; y^2=x(x+2)(x+3);\end{equation}
\begin{equation}\label{ec1}E_3: \; y^2=(x+1)(x+2)(x+3);\quad E_4: \; y^2=x(x+1)(x+2)(x+3);\end{equation}

\noindent $E_0$ and $E_3$ are isomorphic over $\Q$, $E_1,E_2$ are isomorphic over $\Q(i)$, but not over $\Q$, whereas  $E_0,E_3$ are CM curves, and $E_1, E_2,E_4$  are not.
 Therefore,
\begin{equation}\label{p3}n_p(4)=\frac{p}{16}+\frac{a_4(p)}{16} +c_p(4) \;\mbox{for}\; p=4k+3, \end{equation}
\begin{equation}\label{p1}n_p(4)=\frac{p}{16}+\frac{a_0(p)}8+\frac{a_1(p)}8+\frac{a_4(p)}{16}+c_p(4) \;\mbox{for}\; p=4k+1, \end{equation}
where $a_i(p)$ is the (normalized, as before) Frobenius trace of $E_i$; in particular,\linebreak  $a_0(p)=J(k).$
Indeed, for $p=4k+3$ the curves $E_0,E_3$ are {supersingular:}   ${a_0(p)}=0$. On the other hand,
   $E_1$ is a non-trivial quadratic twist of $E_2$ 
which implies  that  $a_1(p)+a_2(p)=0$ for $p=4k+3,$
$a_1(p)=a_2(p) $ for $p=4k+1.$\smallskip

 $\bf{t=5}. $   For $t=5, \;g_5=17$   the same approach yields 15 elliptic curves and a curve of genus $2$, namely:
$$C_{\{1,2,3,4,5\}}:\; y^2=x(x+1)(x+2)(x+3)(x+4).$$
However, counting points on this curve can be reduced to counting points on elliptic curves, since this curve has a non-hyperelliptic involution (see below) and its Jacobian splits (over the base field 
 for $p\equiv 1 \mod 4$ and over its quadratic extension for $p\equiv 3 \mod 4$).    Hence,   $n_p(5)$ can still be expressed in terms of Frobenius  traces for elliptic curves.
\smallskip

 $\bf{t=6}.$ If $t=6 , \;g_6=49$, then some genus 2 curves  $C_T $ for  $  |T|= 5$ have simple Jacobians, which  leads to an expression of  $n_p(6)$ in terms of Frobenius  traces for numerous  ($\ge 35$) elliptic curves and  few ($\le 7$) curves of genus 2.
\smallskip 

  $\bf{t=7}$. For $t=7,\, g_7=129$,  the formulas are of the same type as for $t=6,$ since the only genus 3 curve 
$$C_T, |T|=7, y^2=x(x+1)(x+2)(x+3)(x+4)(x+5)(x+6)$$
has   a non-hyperelliptic involution  and its Jacobian splits into a product of an elliptic curve and the Jacobian of a curve of genus 2. \smallskip 

Note that for $t\le 3$ we have   "explicit" formulas for $n_p(t)$, but, as explained in \cite[section 2.3]{KTVZ}, there is no hope to have similar formulas for $t\ge 4$.

\section{Asymptotic behaviour}

We are now going to show that one can obtain some information on the asymptotic distribution of $n_p(t), p\lra \infty $ from GST, even if no explicit formula for $n_p(t)$ is expected for $t\ge 4$. Let us fix $t\ge1.$

Therefore, the behaviour of $ n_p(t)$ is governed by
the following abelian variety $A$
\begin{equation}\label{at}
A=A(t)=\prod_{T\subset [1,...,t]} J_T, \;J_T=\mathrm{Jac}(C_T) ,\end{equation}
which is the product of all the  Jacobians of the curves
$C_T,$ and thus \\ $\; \dim (A(t))=g_t=2^{t-2}(t-3)+1.$
\medskip

We are interested in the asymptotic distribution (if it exists) of the quantity
$$\delta'_p(t):= n_p(t)-2^{-t}p=\mathrm{trace_A(Fr_p)}+o_p(1) =\sum_{T\subset [1,...,t]}\mathrm{trace_T(Fr_p)}+o_p(1)$$ $ \mathrm{trace_T(Fr_p)}$ being the Frobenius trace on $C_T,$
or rather its normalized version $$\delta_p(t)=\delta'_p(t)/\sqrt{p}\in [-2g_t,2g_t].$$

We begin with 
\subsection{Asymptotics under GS}
As mentioned in section 2,  we suppose now that GST holds for abelian varieties of the form \eqref{at}. Therefore, we get that
$$ \mathrm{ST(A)}\subset\prod_{T\subset [1,...,t]}
\mathrm{ST(J_T)}
$$
Note, however, that this inclusion is always proper for 
any $t\ge 3,$ and thus it is not possible to deduce directly from GST the asymptotics of  the quantity $\delta_p(t).$ 
Nevertheless, for some small $t$ it is possible to refine the analysis and get some equidistribution results  for $\delta_p(t);\;p\lra \infty.$ \medskip

The first non-trivial case is $t=3,\; g_t=1$ which leads to the equidistribution of $\delta_p(3)$ under $\nu_{CM},$
and this is unconditional.\medskip

\subsection{\bf{ Case $t=4$}}
For $t=4,\;g_t=5$ and $A(t)$ takes the form (see \eqref{ec},\eqref{ec1}) 
$$ A(4)=E_0\times E_1\times E_2\times E_3\times E_4\equiv E_0^2\times E_1\times \bar E_1\times E_4$$
 $ \bar E_1$ being the  non-trivial twist of $E_1$ over $\Q(i)$. 
 
 Assuming GST for $B(4):=E_0\times E_1\times E_4$, we get 
 $$ \mathrm{ST(B(4))}\subseteq \Delta_1\times N\times 
\mathrm{SU_2(\C)}\subset \mathrm{\bf{T}}\times \mathrm{SU_2(\C)}\times \mathrm{SU_2(\C)}\subset \mathrm{SU_{6}(\C)},
$$
where $\bf{T}=\mathrm{U_1(\C)}\subset \mathrm{SU_2(\C)},$  $\Delta_1$ being  the image of the diagonal embedding $\bf{T}\hookrightarrow \mathrm{SU_2(\C)}$, whereas $N$ is the
normalizer of  $\bf{T}\hookrightarrow \mathrm{SU_2(\C)}$. In \eqref{p3} and \eqref{p1} the 
summands  $0,0, a_4(p)/16$ and $a_0(p)/8, a_1(p)/8, a_4(p)/16,$ respectively, correspond to the measures

\begin{equation}\label{mu}
\mu_{\Delta_1}= \nu_1(\cdot/2)/2,\;\mu_N=\delta_0+\nu_2/2, \;\nu_{2}. 
\end{equation}\smallskip

\begin{thm}\label{41}
$(i)$ For p running the set $\mathcal P_3:=\{p=4m+3\;\; is \;\ prime\}$ the quantity $\delta_p(4)$ is equidistributed on $[-2, 2]$  with respect to the measure
$$ \mu_3:=\nu_{2}=\frac{\sqrt{4-x^2}dx}{2\pi}.$$

$(ii)$ Assume GST for $B(4)=E_0\times E_1\times E_4.$ 
Then the quantity $\delta_p(4)$ for $p$ running the set $\mathcal P_1:=\{p=4m+1\;\; is \;\ prime\}$ is equidistributed on $[-10, 10]$  with respect to the measure
$$\mu_1:= \lambda_{cm}*\mu*\nu_2 ,$$
where $\mu*\lambda$ denotes the convolution of measures, and 
$$\lambda_{cm}(x)=\frac12 \nu_1(x/2)=\frac{\delta_0}{4}+
\frac{dx}{2\pi \sqrt{16-x^2}}, \mu(x)=\frac12\nu_2(\frac{x}{2})=\frac{\sqrt{16-x^2}dx}{4\pi}.
$$
\end{thm}\smallskip

The proof  follows directly from formulas \eqref{nu},\eqref{mu} and \eqref{p1}.\smallskip

Therefore,

$$\mu_1=  {h_1(x) dx} ,\; h_1=h_{11}+h_{12},\;\hbox{\rm for}$$ 

$$ h_{11}:=\frac{\sqrt{16-s^2}*\sqrt{4-u^2}}{32\pi^2}\ge0,\;h_{12}:=\frac{1}{16\pi^3 \sqrt{16-v^2}}*\sqrt{16-s^2}*\sqrt{4-u^2}\ge0. $$

\smallskip\smallskip

{\em Remark.} Note that the support of both $\lambda_{cm}$ and $\nu_2$ equals the whole interval $[-2,2]$ thus
$$\mathrm{supp}(\lambda_{cm})=\left[- 4, 4\right],\;\; \mathrm{supp}(\mu )= \left[- 4, 4\right],\;\; \mathrm{supp}(\nu_2 )=\left[-2,2\right],$$
and, therefore, $$\mathrm{supp}(\mu_3)=\left[-2, 2\right],\quad \mathrm{supp}(\mu_1)=\left[-10, 10\right] \,.$$
Indeed, if 
$$\mathrm{supp}(f)=[a_1,b_1],\;\;\mathrm{supp}(g)=[a_2,b_2],\;\;\mathrm{supp}(h)=[a_3,b_3]$$
then $\mathrm{supp}(f*g*h)=[a_1+a_2+a_3,b_1+b_2+b_3].$\smallskip\smallskip

{\em Remark.} We have used that $E_0,E_1,E_4$ are not pairwise isogenous.\smallskip

  Clearly,  $E_0$ has CM, and $E_1,\bar  E_1, E_4$ have not, and thus $E_0$   is not isogenous (over $\bar \Q $) to $E_1,\bar  E_1 $ or $ E_4;$  also, $ E_1,\bar  E_1$ are not isogenous to $E_4$ , since $ tr(E'_1),tr(\bar E'_1)=\pm 2$ and $tr E_4'=-1$ for $E'_k:=E_k \mod 5,$ and the traces are not equal up to sign.  
 \qed\smallskip

Note that  GST($B(4)$)  implies that the sequences 
$\{a_0(p)\},\{a_1(p)\} $ and $\{a_4(p)\}$ are statistically independent.

  The condition in Theorem {\rm\ref{41}},{\rm (ii)} is quite plausible, as mentioned in Section 3.

\smallskip

Looking at the support of $\mu_1$ and $\mu_3$, we obtain

\begin{cor}\label{c} Under the conditions of Theorem {\rm\ref{41} } for any $\eps>0$  there exist $4$ primes
$p_1=p_1(\eps)\equiv 1\mod 4$, $p_3=p_3(\eps)\equiv 3\mod 4,$  
$p'_1=p_1'(\eps)\equiv 1\mod 4$, $p'_3=p'_3(\eps)\equiv 3\mod 4,$
such that 
$$n_{p_1}(4)\ge \frac{p_1}{16}+\left(\frac58-\eps\right) \sqrt{p_1},\;\;n_{p'_1}(4)\le \frac{p'_1}{16}-\left(\frac58-\eps\right) \sqrt{p'_1} , $$
$$n_{p_3}(4)\ge \frac{p_3}{16}+\left(\frac18-\eps\right) \sqrt{p_3},\;\;n_{p'_3}(4)\ge \frac{p'_3}{16}-\left( \frac18 -\eps\right) \sqrt{p'_3}. $$
    \end{cor}

This corollary shows that the upper bounds for $ n_p(4)-p/16$ are essentially tight. Note that the existence of 
$p_3$ and $p'_3$ is unconditional.

\subsection{\bf Case $t=5$}
For the next value $t=5$ one has $g(5)=17$ and $A(5)$ is (isogenous to)
a product of 17 ECs; however, many of them are isomorphic or conjugate.

Indeed,

$$ A(5)=A(4)\times E_5\times E_6\times E_7\ldots\times E_{14}\times  \mathrm{Jac}(C)$$
for
$$ A(4)=E_0\times E_1\times E_2\times E_3\times E_4\equiv E_0^2\times E_1\times \bar E_1\times E_4$$

\noindent with  $E_0,E_1,E_2,E_3$ and $E_4$ defined above;

$$E_5: y^2=x(x+1)(x+4);\; E_6: y^2=x(x+2)(x+4);\;E_7: y^2=x(x+3)(x+4);$$
$$E_8: y^2= (x+1)(x+2)(x+4);\, E_9: y^2=(x+1)(x+3)(x+4);\,E_{10}: y^2=(x+2)(x+3)(x+4); \,$$ 
$$ E_{11}: y^2=x(x+1)(x+2) (x+4);\,\;\; E_{12}: y^2=x(x+1) (x+3)(x+4);\; $$ 
$$E_{13}: y^2=x(x+2)(x+3)(x+4);\;\; E_{14}: y^2=(x+1)(x+2)(x+3)(x+4),\; $$
  $$C:  y^2=F(x)=x(x+1)(x+2)(x+3)(x+4),$$
  
\noindent note that the curve $C$ admits a non-hyper-elliptic involution 

$$\sigma:(x,y)\mapsto \left(-\frac{2x+6}{x+2},\;\frac{2\sqrt 2 y}{(x+2)^3}\right) $$

\noindent since $(x+2)^6F\big(\frac{2}{x+2}\big)=8F(x),$
and thus its trace can be calculated from the traces of a pair 
of elliptic curves, say, $E_{15}$ and $E_{16}$, where

$$E_{15}:y^2=x^3+2\sqrt 2x^2-9x-18\sqrt 2,$$

$$E_{16}:y^2=x^3+2 x^2-9x/2 - 9,$$ 

\noindent $E_{15},E_{16}$ being conjugate over $\Q(\sqrt 2).$ \smallskip

One also notes numerous isomorphisms between those curves. \smallskip

Namely,

 \smallskip

\noindent $E_0$ is $\Q$-isomorphic to $E_3,\;{ {E_6}'},\;E_{10}; j(E_0)=1728=2^63^3;$\smallskip\smallskip

\noindent $ E_{4} $ is $\Q$-isomorphic to $ E_{5},\; \bar E_{7},\;       E_{14},   \;j(E_4)=35152/9=2^43^{-2}13^3;$  \smallskip\smallskip

\noindent $E_{1}$ is $\Q$-isomorphic to  $\Bar{ E_{2}}, \;  E_8,\;\Bar E_9 \; E_{11}, \;E_{13};\;j(E_1)=21952/9=2^63^{-2}7^3; $ \smallskip\smallskip

 \noindent $J(E_{12}) = 1556068/81=2^23^{-4}73^3; $\smallskip\smallskip

\noindent$E_{15}\simeq_{\Q} \;   {E_{16}'},j(E_{15})=2744000/9=2^63^{-2}5^37^3,$ \smallskip\smallskip

\noindent where ${\Bar {E}}$ is the $ \Q(i)$-twist of $ E,$ and ${ {E}'}$ is the $\Q(\sqrt{2})$-twist of $E.$

Therefore,
$$ A(5)=E_0^3\times   E_0'\times E_1^4\times \bar E_1^2 \times  E_4^3\times \bar E_{4}\times E_{12}\times E_{15}\times {E_{15}'},$$
where all elliptic curves in the  product are pairwise non-isogenous.

\smallskip Indeed, $E_0,E_1$ and $E_4$ are pairly non-isogenous (over $\bar \Q$). For the pair $(E_1,E_{12})$ it is sufficient to look at  their reductions modulo 7 to see that $E''_1=E_1 \mod 7$ is supersingular, $a_1(7)=0,$ whereas $E''_{12}=E_{12} \mod 7$   is not, $a_{12}(7)=-1.$ Thus $E_1''$ and $E''_{12}$ are not isogenous over $\bar \F_7$ and hence $E_1, E_{12}$ are not isogenous over $\bar \Q $. The same for $(E_1,E_{16})$ since $a_{16}(7)=5.$ For $(E_4,E_{12}): a_{4}(11)=5, 
 a_{12}(11)=-3;$ for $(E_4,E_{16}): a_{4}(7)=1, 
 a_{16}(7)=5;$ finally, $a_{12}(7)=-1,a_{16}(7)=5 $ gives the same for 
 $(E_{12},E_{16})$.
 \smallskip \smallskip

Admitting GST for $A(5),$ one gets a result parallel to Theorem {\rm\ref{41} } where the formulas are more lengthy and include  more convolutions; in fact it is sufficient to admit  GST for
$$ B(5):=E_0\times   E_1\times    E_4\times    E_{12}\times E_{16}, \;\dim B(5)= 5. $$

The results depend on $p$ modulo 8; for instance, in the simplest case $p=8m+7$ one has
$$ n_p(5)-\frac{p}{32}=\frac{a_0(p)}{16}+ \frac{a_1(p)}{16}+\frac{a_4(p)}{16}+\frac{a_{12}(p)}{32}, $$\smallskip

\noindent since $p$ is inert both in $\Q(\sqrt{-1})$ and  $\Q(\sqrt{2})$ which gives the simplification of all other  terms (as for $p=4k+3$ if $t=4$).

\smallskip In general we get, using how $p$ splits in  $\Q(\sqrt{-1})$ and  $\Q(\sqrt{2})$:

\begin{thm}\label{42}

$(i)$  Admit GST for
$$ B'(5):=E_0\times   E_1\times    E_4\times    E_{12} , \;\dim B'(5)= 4. $$

If $p=8m+7, $ then

$$\big|\delta(p)\big|\le \frac{7}{32}\;\;\;  \hbox{ \rm{for}}\;\;\; \delta(p):=\frac{1}{2\sqrt p}\Big| n_p(5)-\frac{p}{32}\Big|.$$

 On the sequence $\mathcal {P}_7=\{p=8m+7\}$ the quantity $\delta_7(p)$ is equidistributed for the measure
$$\mu_7(5)=\lambda_{cm}*\mu*\mu*\nu_2 $$ 
in the notation of Theorem {\rm\ref{41} }.\smallskip\smallskip

$(ii)$  Admit GST for $B'(5).$

If $p=8m+5, $ then 
$$ n_p(5)-\frac{p}{32}=\frac{a_0(p)}{16}+ \frac{a_1(p)}{8}+\frac{3a_4(p)}{16}+\frac{a_{12}(p)}{32}, $$ 
and
$$\big|\delta(p)\big|\le \frac{13}{32}\;\;\;  \hbox{ \rm{for}}\;\;\; \delta(p):=\frac{1}{2\sqrt p}\Big| n_p(5)-\frac{p}{32}\Big|.$$

On the sequence $\mathcal {P}_5=\{p=8m+5\}$ the quantity $\delta(p)$ is equidistributed for the measure
$$\mu_5(5)=\lambda_{cm}*\mu*\mu'*\nu_2 $$
for $\mu':=\mu(x/3)/3. $\smallskip\smallskip

$(iii)$   Admit GST for $B(5).$

If $p=8m+3, $ then 
$$ n_p(5)-\frac{p}{32}=\frac{a_0(p)}{8}+ \frac{a_1(p)}{16}+\frac{a_4(p)}{16}+\frac{a_{12}(p)}{32}+\frac{a_{15}(p)}{16}, $$ 
and
$$\big|\delta(p)\big|\le \frac{11}{32}\;\;\;  \hbox{ \rm{for}}\;\;\; \delta(p):=\frac{1}{2\sqrt p}\Big| n_p(5)-\frac{p}{32}\Big|.$$

On the sequence $\mathcal {P'}_3=\{p=8m+3\}$ the quantity $\delta(p)$ is equidistributed for the measure
$$\mu_3(5)=\lambda'_{cm}*\mu*\mu*\nu_2*\mu $$
for $\lambda'_{cm}:=\lambda_{cm}(x/2)/2. $\smallskip\smallskip

$(iv)$   Admit GST for $B(5).$

If $p=8m+1, $ then 
$$ n_p(5)-\frac{p}{32}=\frac{a_0(p)}{8}+ \frac{a_1(p)}{8}+\frac{3a_4(p)}{16}+\frac{a_{12}(p)}{32}+\frac{a_{15}(p)}{16}, $$ 
and
$$\big|\delta(p)\big|\le \frac{17}{32}\;\;\;  \hbox{ \rm{for}}\;\;\; \delta(p):=\frac{1}{2\sqrt p}\Big| n_p(5)-\frac{p}{32}\Big|.$$

On the sequence $\mathcal {P'}_1=\{p=8m+1\}$ the quantity $\delta(p)$ is equidistributed for the measure
$$\mu_1(5)=\lambda'_{cm}*\mu*\mu'*\nu_2*\mu. $$
\end{thm}

\begin{cor} In the conditions of Theorem {\rm\ref{42} } .\smallskip

$(i)$ One has  the bound
 $$\big|\delta_7(p)\big|\le \frac{17}{32}   $$ 
which is tight on $\mathcal {P'}_1$.\smallskip

$(ii)$  On $\mathcal {P'}_3$ one has
$$\big|\delta_7(p)\big|\le \frac{13}{32} ,  $$ 
tight on $\mathcal {P'}_3$.\smallskip

$(iii)$  On $\mathcal {P}_5$ one has
$$\big|\delta_7(p)\big|\le \frac{11}{32} ,  $$ 
tight on $\mathcal {P}_5$.\smallskip

$(iv)$  On $\mathcal {P}_7$ one has
$$\big|\delta_7(p)\big|\le \frac{7}{32} ,  $$ 
tight on $\mathcal {P}_7$.
\end{cor}

\section{Unconditional results}

 \subsection{\bf $t$=4.}   Note than that in \cite{Jo} the GST(A) is proven, in particular, for any product $A=E\times F$  of two (non-isogenous) elliptic curves. Using this result, one can obtain an unconditional result similar to (but weaker than) Corollary \ref{c}. Indeed, one has 

\begin{lemma} The sequences
$$S_0=\{a_0(p)\},\; S_1=\{a_1(p)\},\;S_4=\{a_4(p)\} $$
are pairwise  statistically independent.
\end{lemma}

Indeed, the GST holds for $A_{01}=E_0\times E_1, A_{04}=E_0\times E_4,$ and $A_{14}=E_1\times E_4$ thanks to \cite{Jo}, which immediately implies that independence.

\begin{prop}\label{3}   For any $\eps>0$  there exist  primes $p_1,p_1' \,$ with
$$p_1=p_1(\eps)\equiv 1\mod 4,\; p_1=p_1'(\eps)\equiv 1\mod 4,$$  
\noindent such that 
$$n_{p_1}(4)\ge \frac{p_1}{16}+\left(\frac38-\eps\right) \sqrt{p_1},\;\;n_{p'_1}(4)\le \frac{p'_1}{16}-\left(\frac38-\eps\right) \sqrt{p'_1} . $$
\end{prop}

\emph{Proof.} Fix $\eps>0.$ Since $S_0$ and $S_1$ are independent,
one can find $p_1,p_1'$ s.t. 
$$\frac{a_0(p_1)}{\sqrt{p_1}}\ge 2-4\varepsilon,\;\frac{a_1(p_1)}{\sqrt{p_1}}\ge 2-4\varepsilon,$$ 
$$\frac{a_0(p'_1)}{\sqrt{p_1}}\le 4\varepsilon-2,\;\frac{a_1(p_1)}{\sqrt{p'_1}}\le 4\varepsilon-2,$$ 
\noindent and since $|{a_4(p_1)}|\le 2\sqrt{p_1}$, an immediate calculation implies the conclusion.\qed\smallskip\smallskip

 \subsection{\bf $t$=5.} The same GST$(E\times F)$    can also be applied to get  unconditional results for $t=5$ in the spirit of Proposition \ref{3} which are rather  weak, yet non-trivial. 

 \begin{prop}\label{7} 
$(i)$ For any $\eps>0$  there exist  primes $p_1,p_1' \,$ with
$$p_1=p_1(\eps)\equiv 1 \mod 8,\; p_1=p_1'(\eps)\equiv 1 \mod 8,$$ 

\noindent such that  

$$n_{p_1}(5)\ge \frac{p_1}{32}+\left(\frac{1}{16}-\eps\right) \sqrt{p_1},\;\;n_{p'_1}(5)\le \frac{p'_1}{32}-\left(\frac{1}{16}-\eps\right) \sqrt{p'_1} . $$
 
 $(ii)$ For any $\eps>0$  there exist  primes $p_3,p_3' \,$ with
$$p_3=p_3(\eps)\equiv 3\mod 8,\;p'_3=p_3'(\eps)\equiv 3\mod 8,$$ 

\noindent such that  

$$n_{p_3}(5)\ge \frac{p_3}{32}+\left(\frac{3}{16}-\eps\right) \sqrt{p_3},\;\;n_{p'_3}(5)\le \frac{p'_3}{32}-\left(\frac{3}{16}-\eps\right) \sqrt{p'_3} . $$
 
 $(iii)$ For any $\eps>0$  there exist  primes $p_5,p_5' \,$ with
$$p_5=p_5(\eps)\equiv 5\mod 8,\;p'_5=p_5'(\eps)\equiv 5\mod 8,$$ 

\noindent such that  

$$n_{p_5}(5)\ge \frac{p_5}{32}+\left(\frac{1}{16}-\eps\right) \sqrt{p_5},\;\;n_{p'_5}(5)\le \frac{p'_5}{32}-\left(\frac{1}{16}-\eps\right) \sqrt{p'_5} . $$\pagebreak
 
 $(iv)$ For any $\eps>0$  there exist  primes $p_7,p_7' \,$ with
$$p_7=p_7(\eps)\equiv 7\mod 8,\;p'_7=p_7'(\eps)\equiv 7\mod 8,$$ 

\noindent such that  

$$n_{p_7}(5)\ge \frac{p_7}{32}+\left(\frac{1}{16}-\eps\right) \sqrt{p_7},\;\;n_{p'_7}(5)\le \frac{p'_7}{32}-\left(\frac{1}{16}-\eps\right) \sqrt{p'_7} . $$
\end{prop}

Indeed, e.g., in $(iv)$ 
$$ n_p(5)-\frac{p}{32}=\frac{a_0(p)}{16}+ \frac{a_1(p)}{16}+\frac{a_4(p)}{16}+\frac{a_{12}(p)}{32}, $$

\noindent and one can choose  first two terms very close to $\sqrt{p}/8$  or to $- \sqrt{p}/8,$ by independence of ${a_0(p)}$  and  ${a_1(p)}$. Similar arguments work for all other cases.

\end{document}